# ON THE COLLATZ PROBLEM

by
Roupam Ghosh
roupam_ghosh@vsnl.net

**Abstract:**
Taking a new approach towards analyzing the Collatz Problem, or, 3x+1 conjecture. Introducing some new functions, the Collatz-2 and Collatz-3 sequences, as well as deducing results related to Collatz-2 and Collatz-3 sequences.

## The Collatz Problem

Consider the function

$$f(x) = x/2 \quad \text{if } x \text{ is even}$$
$$\quad = 3x+1 \quad \text{if } x \text{ is odd}$$

The Collatz conjecture states that there exists a number ‘d’ corresponding to x such that $f^d(x) = 1$ for all natural numbers x

where $f^d(x) = f(f(f(\ldots \text{ d times } \ldots(f(f(f(x)))))))$

The sequence of numbers following the collatz function is as:

| | |
|---|---|
| for x=1 | 1 |
| for x=2 | 2 1 |
| for x=3 | 3 10 5 16 8 4 2 1 |
| for x=4 | 4 2 1 |
| for x=5 | 5 16 8 4 2 1 |
| for x=6 | 6 3 10 5 16 8 4 2 1 |

...

## The Collatz-2 sequence

Now we construct a new sequence:
We start normally from 1 as in Collatz Sequence defined above
If we reach a number that already exists in any sequence before, then we stop iterating
Then we have, the **collatz-2** sequence, as follows

| | | |
|---|---|---|
| for x=1 | 1 | |
| for x=2 | 2 | (Since 1 already exists before) |
| for x=3 | 3 10 5 16 8 4 | (Since 2 already exists before) |
| for x=4 | | (Since 4 already exists before) |
| for x=5 | | (Since 5 already exists before) |
| for x=6 | 6 | (Since 3 already exists before) |
| for x=7 | 7 22 11 34 17 52 26 13 40 20 | (Since 10 already exists before) |

...

Here each row is called a level

This collatz-2 sequence brings computability into question
**We can never know about a particular level unless we compute all levels above it**

We now make a couple of definitions based upon the above sequence

| | | |
|---|---|---|
| touch(N) | = | No. of numbers covered before N in the collatz-2 sequence |
| level(N) | = | x if xth series in collatz-2 sequence contains N |
| s(N) | = | Steps required to reach 1 starting from N |
| e(N) | = | Number of elements in the Nth level |
| max(N) | = | The maximum number in the Nth level in Collatz-2 sequence |

For example in the collatz-2 sequence

| | | |
|---|---|---|
| touch(4) | = | 7 |
| level(4) | = | 3 |
| s(4) | = | 3 |
| e(4) | = | 0 |

We shall call the function e(N) as epsilon of N

The first number and last number in any level will be called level starters and level enders repectively for that level.

We shall study this Collatz-2 sequence.

In general we assume the following to be true:

1. If a non-trivial cycle exists in the Collatz-2 sequence, then it also exists in the Collatz sequence
2. Also, if any level does not stop iterating in the Collatz-2 sequence, then it also does not stop iterating in the Collatz sequence.

## The Lambda Function

**Consider** the lambda function as follows

$$\Lambda(m) = 0 \text{ if } e(m) = 0$$

$$\Lambda(m) = 1 \text{ if } e(m) > 0$$

then, we call nz(n) the number of non-zero lambdas for numbers less than or equal to n

$$nz(n) = \sum \Lambda(i) \qquad \text{for } i = 1 \text{ to } n$$

Let z(n) represent the number of zeroes of lambdas for number less than or equal to n

then $z(n) + nz(n) = n$

Now,

$nz(n)/n = (1/n) * \sum \Lambda(i)$ for $i = 1$ to $n$

$\lim nz(n)/n < 1$ as $n \to \infty$

Now let us observe something,
If we consider the collatz-2 sequence then we have the following:

## Some Lemmas:

**Lemma 1: If level(n) < n then e(n) = 0**

If $e(n) = 0$ then definitely there exists one $m < n$ such that $f^d(m) = n$ and $d < e(m)$
Hence $\text{level}(n) = m < n$

**Lemma 2: If e(n) ≠ 0 then level(n) = n**

If $e(n) \neq 0$ then definitely n has at least one member in its level
Since level(n) starts with n, hence $\text{level}(n) = n$

Combining both we get
**level(n) ≤ n for all natural numbers n**

**Lemma 3: $e(x) = s(x) - s(f^{e(x)}(x))$**
This follows from the basic contruction of the Collatz-2 sequence
In any level
if $e(x) = 0$ then $s(x) - s(f^0(x)) = s(x) - s(x) = 0$
if $e(x) > 0$ then
$s(x) - s(f^{e(x)}(x))$ = the number of elements in level(x) = e(x)

**Lemma 4: $s(x) > s(f(x)) > s(f^2(x)) > ... > s(f^d(x)) > ... > 0$**
Since x comes before f(x) in the collatz sequence, hence $s(x) > s(f(x))$

**Lemma 5: $e(f^n(x)) = 0$ if level(x) = x and n < e(x)**

Consider level(x)
Then obviously, it looks like

$x, f(x), f^2(x),... , f^{e(x) - 1}(x)$
For $1 \leq n < e(x)$, we have $e(f^n(x)) = 0$
Hence, proved

**Lemma 6: $s(2^n k) = s(k) + n$**
For,
$s(2^n k) = s(2^{n-1}k) + 1 = ... = s(k) + n$

**Lemma 7:** **$\sum e(n) < \max(m_1,m_2,m_3,...,m_n)$ where $m_i = \max(i)$, $level(m_i) \leq n$, and $e(n) > 1$**

We shall call $m_n$ as maximal(n)
Now, it is obvious, that $\sum e(n) \leq z(maximal(n))$
and, $z(maximal(n)) < maximal(n)$ for $e(n) > 1$
Hence, follows

**Lemma 8:** **$(1/n)\sum \Lambda(i).level(i) < (n+1)/2$**

i.e., Average of level starters is less than average of n

**Lemma 9:** **$\lim nz(n)/n < 0.5$ as $n \to \infty$ if Collatz conjecture is true**

If collatz problem is true

then as $n \to \infty$ $nz(n) < z(n)$

ie., $nz(n) < n - nz(n)$

ie., $nz(n)/n < 1/2$

# Cycles and their interpretation in the Collatz-2 sequence

Now lets consider 2 cases:
1. For all x, e(x) is finite
2. For a particular x=n, e(n) is infinite and e(x) is finite for all x < n

## Case 1

Since e(x) is finite, includes 2 possibilities,
a. $level(f^{e(x)}(x)) < x$
b. $level(f^{e(x)}(x)) = x$

Here possibility a. denotes a normal level, which satisfies the Collatz problem
Possibility b. on the other hand denotes the existence of a cycle which contradicts the Collatz problem

## Case 1, b

If we have cycles in a particular level(n), then Collatz-2 sequence for that level stops iterating after a finite time. This comes from the fundamental property of the Collatz-2 sequence ie., "No number is repeated"

Now lets analyze case 1, b in details.
Since a cycle exists at a particular level, lets say the level starter is x
and for some $1 < k < e(x)$
$f^k(x) = f^{e(x)}(x)$

Lets now consider the following possibilities, where $f^k(x) = f^{e(x)}(x)$

| | $f^{k-1}(x)$ | $f^k(x)$ | $f^{e(x)-1}(x)$ | $f^{e(x)}(x)$ |
|---|---|---|---|---|
| 1 | odd | even | even | even |
| 2 | odd | even | odd | even |
| 3 | even | even | even | even |
| 4 | even | even | odd | even |
| 5 | even | odd | even | odd |

Other possibilities do not exist.
Note that cycleing starts from $f^{e(x)}(x)$, where $f^{e(x)}(x) = f^k(x)$

Consider possibility 1:
Let, $f^{k-1}(x) = d$ where d is odd
then $f^k(x) = 3d+1$

Let, $f^{e(x)-1} = e$ where e is even
then, $f^{e(x)} = e/2$

Here, $e/2=3d+1$

Consider possibility 2:
Let, $f^{k-1}(x) = d$ where d is odd
then $f^k(x) = 3d+1$

Let, $f^{e(x)-1} = e$ where e is odd
then, $f^{e(x)} = 3e+1$

Here, $3d+1 = 3e+1$
Hence, $d=e$

But that is not possible since cycling starts from $f^{e(x)}$ and not $f^{e(x)-1}$
Hence we cancel possibility 2.

Consider possibility 3:
Let, $f^{k-1}(x) = d$ where d is even
then $f^k(x) = d/2$

Let, $f^{e(x)-1} = e$ where e is even
then, $f^{e(x)} = e/2$

Here, $d/2 = e/2 \Rightarrow d = e$
Hence we have the same case as possiblity 2
So we cancel out possibility 3

Consider possibility 4:
Let, $f^{k-1}(x) = d$ where d is even

then $f^k(x) = d/2$

Let, $f^{e(x)-1} = e$ where e is odd
then, $f^{e(x)} = 3e+1$

$d/2 = 3e+1$

Consider possibility 5:
Let, $f^{k-1}(x) = d$ where d is even
then $f^k(x) = d/2$

Let, $f^{e(x)-1} = e$ where e is even
then, $f^{e(x)} = e/2$

Here, also $d/2 = e/2$
Hence, $d = e$

So we cancel out possibility 5

We are left with possibility 1 and possibility 4

In both cases, we can see that,
An even number and an odd number generates the same number in the same level, ie, x

Hence **if we can show that an even number and an odd number cannot generate the same number in the same leve**l, we prove that cycles don't exist in the Collatz-2 sequence, and hence the Collatz sequence.

We call such an even and odd number pair as the **"cyclic pair"**

**Interpretation of cyclic pairs:**

Consider,
In level x, o and e are odd and even numbers
and they generate the same number,
i.e,

$f(o) = f(e)$

ie., $3o+1 = e/2$

We have two cases
1. $o > e$
2. $e > o$

Case 1.

if $o > e$
Then,
$3o + 1 > e/2$

But here, $3o+1 = e/2$

Hence case 1. is not valid

Case 2.

if e > o
Since e > o hence, e is definitely not the level starter
We have the following subcases

Subcase 1.

*o is the level starter, e is the level ender*

Subcase 2.

*o is not the level starter, e is the level ender*
*The level starter must be an odd integer,*
*since the level accommodates more than*
*two integers.*

Subcase 3.

*e is not the level starter, o is the level ender*
*The level starter must be an odd integer,*
*since the level accommodates more than*
*two integers.*

To analyze Case 2. we shrink the Collatz-2 sequence into Collatz-3 sequence
The rules for the Collatz-3 sequence is as follows

Remove all even integers from Collatz-2 sequence to get the Collatz-3 sequence.

The sequence looks as follows
**Collatz-2 sequence**

for x=1 1
for x=2 2
for x=3 3 10 5 16 8 4
for x=4
for x=5
for x=6 6
for x=7 7 22 11 34 17 52 26 13 40 20
...

**Collatz-3 sequence**

for x=1 1
for x=2
for x=3 3 5
for x=4
for x=5
for x=6
for x=7 7 11 17 13
...

It is obvious, that if there exists a cycle in the Collatz-2 sequence, then it also exists in the Collatz-3 sequence

Let $a_1$ $a_2$ ... $a_m$ denote an endless cycle, from $a_1$ to $a_m$ and back to $a_1$
Let the odd numbers in the cycle be $o_1$ $o_2$ $o_3$ ... $o_k$ where $k < m$
Then the Collatz-3 sequence will have the same number of odd numbers cycling,
Because Collatz-3 applies the function

$f(x) = (3x+1)/2^k$ when x is odd,
and k is the maximum power of 2 that divides 3x+1

So stripping Collatz-2 of even numbers doesn't matter with regards to cycle,
Since,
$o_k = (3o_{k-1} + 1)/2^r$ for some r

$a_m = o_k$ if $a_m$ is odd
$= (3o_k+1)/2^p$ if $a_m$ is even, for some p

$o_1 = (3o_k+1)/2^q$ for some $q > p$ Since, $a_m$ cycles to $a_1$

Now we again consider the three subcases:

Subcase 1.

o is the level starter, e is the level ender
which in Collatz-3 sequence will be

*$o = o_1$ is the level starter and $o_k$ is the level ender*
*such that $(3o_k + 1)/2^r = e$ for some r*

*Here cycle is from $o_1$ to $o_k$ and back to $o_1$*

Subcase 2.

*$o_1$ is the level starter and $o_k$ is the level ender*
*such that $(3o_k + 1)/2^r = e$ for some r*
*and $o_1 \neq o$, hence $o > o_1$*
*$o = o_j$ for some j, $1 < j < k$*
*Here cycle is from $o_1$ to $o_k$ and back to $o_j$*

Subcase 3.

*$o_1$ is the level starter and $o = o_k$ is the level ender*
*and $o_1 \neq o$, hence $o > o_1$*
*Here cycle is from $o_1$ to $o_k$ and back to $o_j$*
*for some j, $1 < j < k$*

In subcase 1, the cycle is from end to the begining of the level

In subcase 2 and 3, the cycle is from end to somewhere between the begining and end.

So, if we prove that the cycle cannot occur among integers (odd) in the Collatz-3 sequence, then it implies that cycles cannot occur in the Collatz-2 and hence Collatz sequence.

**Analysis of the endless cycling in Collatz-3 sequence**

Consider the sequence, $o_1, o_2, \ldots , o_n$ where each element is odd

If this is a cycle, then

$o_2 = (3o_1+1)/2^{k(1)}$

$o_3 = (3o_2+1)/2^{k(2)}$

...

$o_n = (3o_{n-1}+1)/2^{k(n-1)}$

$o_{n+1} = (3o_n+1)/2^{k(n)}$

But, $o_{n+1} = o_1$

Hence,

$o_1 = (3o_n+1)/2^{k(n)}$

Therefore we get

$k(1) = (1/\log2)(\log(3o_1+1) - \log(o_2))$

$k(2) = (1/\log2)(\log(3o_2+1) - \log(o_3))$

...

$k(n) = (1/\log2)(\log(3o_n+1) - \log(o_1))$

Hence

$\sum k(n) = (1/\log2)\sum(\log(3o_n+1) - \log(o_n+1))$

$> (1/\log2)\sum(\log(3o_n) - \log(o_n+1))$

$= n (\log3/\log2) + (1/\log2) \sum(\log(o_n) - \log(o_n+1))$

$= n (\log3/\log2) + 0$ [ Since $o_{n+1} = o_1$ ]

Hence,

$\sum k(n) > n(\log3/\log2)$, if there is a cycle

Here, $\sum k(n)$ is the sum of powers of 2.

Hence, in the corresponding Collatz-2 sequence, it represents the number of “divide by 2” operations in the repective level, and n is the number of “mutiply by 3 and add 1” operations till $o_n$

Hence,

$\sum k(n)/n > \log3/\log2$

But,

$\sum k(n)+n = e(m)$ where m is the corresponding level

Therefore,if we denote number of “divide by 2” operations in the cycle as q and if we denote number of “multiply by three and add 1” operations in the cycle as p

Then,

$p < e(m)/(1+\log3/\log2)$
$q > e(m)/(1+\log2/\log3)$
$p+q = e(m)$

## Case 2:

If e(n) is infinite, then we can never know about e(n+1) until we reach n+1 in level(n), in which case e(n+1) = 0

In a similar way, for any k > n, we cannot know about e(k) unless we reach k in level(n), in which case e(k) = 0

So if e(n) is infinite, then our information is limited and based upon computational evidence. We can never complete computing for level(n), hence the collatz-2 sequence will have no information for k > max(level(n))

Hence, our knowledge of collatz-2 sequence will increase with time, but we can never complete constructing the sequence if such a level exists where e(x) is infinite.

In this case, we can say the Collatz-2 sequence is incomplete.
Thus, analysis of the Collatz-2 sequence fails, and the theorem remains a conjecture.

# Distribution of zeroes of lambda function

We denote
z(n) as the number of zeroes of the lambda function for $x \leq n$
nz(n) as the number of ones of the lambda function for $x \leq n$

z(n) and nz(n) are monotonically increasing functions

Now, $\sum e(n) > n$ then since $\sum e(n)$ is a monotonically increasing function, the probability of the truth of Collatz conjecture increases with increasing n.
In any case,

$\lim \sum e(n)/n > 1$ as $n \to \infty$

Consider,

$\sum e(n) - nz(n) = z(n)$ + number of zeroes after n and less than maximal(n)

Now, number of zeroes after n and less than maximal(n) $\leq$ z(maximal(n)) - z(n)

Hence,

$\sum e(n) - nz(n) \leq z(maximal(n))$, whenever $e(n) > 1$

But we already know, that for $e(n) > 1$

$z(maximal(n)) \geq \sum e(n)$

For sufficiently large n,

$$z(\text{maximal}(n)) \sim \sum e(n)$$

ie., $nz(n)/n \to 0$ as $n \to \infty$

Hence, for sufficiently large values of n, $z(n) \sim n$